\documentclass[12pt]{article}
\usepackage[ansinew]{inputenc}
\usepackage{pst-all,epsfig}
\usepackage{amsmath}
\usepackage{amsfonts}
\usepackage{float}
\usepackage{graphicx, graphics}
\usepackage{amssymb}
\usepackage[left=2.5cm,top=2cm,right=2.5cm,bottom=2.5cm]{geometry}
\newcommand{\fin}{\hfill $\Box$}

\usepackage[USenglish]{babel}

\usepackage{graphics}
\newcommand{\R}{\mathbb{R}}

\newcommand{\Rt}{\mathbb{R}^3}
\newcommand{\Rn}{\mathbb{R}^{n}}

\newcommand{\Sn}{\mathbb{S}^{n-1}}
\newcommand{\inte}{\operatorname*{int}}

\newcommand{\aff}{\operatorname*{aff}}

\newcommand{\bd}{\operatorname*{bd}}

\newcommand{\lin}{\operatorname*{lin}}

\newcommand{\diam}{\operatorname*{diam}}

\newcommand{\GL}{\operatorname*{GL}}

\newtheorem{lemma}{Lemma}
\newtheorem{theorem}{Theorem}

\title{A characterization of centrally symmetric convex bodies in terms of visual cones}
\author{
E. Morales-Amaya$^{1}$, J. Jer\'onimo-Castro$^{2}$, \\ and D. J. Verdusco Hern\'andez$^{3}$ \\ \\
\small{$^{2}$Facultad de Ingenier\'ia}\\
\small{Universidad Aut\'onoma de Quer\'etaro, M\'exico}\\
\small{$^{1-3}$Facultad de Matem\'aticas-Acapulco,}\\
\small{Universidad Aut\'onoma de Guerrero, M\'exico}\\
 \small{\texttt{$^{1}$emoralesamaya@gmail.com,$^{2}$jesus.jeronimo@uaq.mx,}}\\
 \small{\texttt{$^{3}$diana.janett.h@gmail.com}}
}

\date{\small{\today}}

\begin{document}
 \maketitle
\begin{abstract}  
In this work we prove the following result: Let $K$ be a strictly convex body in the Euclidean space $\mathbb{R}^n, n\geq 3$,  and let $L$ be a hypersurface, which is the image of an embedding of the sphere $\mathbb{S}^{n-1}$, such that $K$ is contained in the interior of $L$. Suppose that, for every $x\in L$, there exists $y\in L$ such that the support double-cones of $K$ with apexes at $x$ and $y$, differ by a translation. Then $K$ and $L$ are centrally symmetric and concentric.
\end{abstract}
 
\section{Introduction}
Let $K\subset \mathbb{R}^{n}$ be a convex body, i.e., a compact and convex set with non-empty interior, $n \geq 3$, and let $x\in \mathbb{R}^{n}\setminus K$. We call the set 
 \[
\bigcup _{y\in  K} \aff \{x,y \}
\]
the \textit{solid double-cone} generated by $K$ and $x$, where $\aff\{x,y\}$ denotes the affine hull of $x$ and $y$. The boundary of the solid double-cone generated by $K$ and $x$ will be called the support double-cone of $K$ with apex at $x$ and it will be denoted by $C_x$.  In what follows, we shall use the names cone or support cone, by simplicity, to refer to the double cone and double support cone, respectively.

A classical problem in convexity is to determine properties of a convex body 
$K\subset \Rn$ from the information of its orthogonal projections. For instance, in dimension 3 one can prove that if all the orthogonal projections of a body $K$  are circles, then $K$ is a Euclidean ball. One can see this problem from the following perspective: consider the family of cylinders where $K$ is inscribed and impose a condition in a particular section of each of them, which is obtained with a hyperplane perpendicular to the lines which generates the cylinder.   
In our example, this means that we have a convex body $K\subset \Rt$ such that for every cylinder $\Omega$, where $K$ is inscribed, the section $H\cap \Omega$ is a circle, where $H$ is a plane orthogonal to the lines which determines $\Omega$. 

We formulate the following general problem.

\textbf{Problem 1.} \emph{For a given subgroup $G$ of the general linear group $\GL(\mathbb{R},n)$, to determine the convex bodies $K\subset \Rn$, $n\geq 3$,  such that for every couple of different cylinders $\Lambda, \Gamma$, circumscribed to $K$, there exists an element $\Phi \in G$ such that $\Phi(\Lambda)= \Gamma$}. 

Kuzminyh \cite{Kuzminyh} proved, for $n=3$, that the assumption $G=O(\R,3)$ implies that $K$ is a ball, where $O(\R,3)$ is the real orthogonal group. On the other hand, if $K\subset \Rn, n\geq 3$, is centrally symmetric, in virtue of the Aleksandrov Uniqueness Theorem, it follows that $K$ is a ball since all the projections have the same volume. Recently, L. Montejano \cite{Montejano} has considered the case where $G$ is the affine subgroup $A(\R,n)$ and he has obtained that $K$ is an ellipsoid.

In virtue that the cylinders are cones with apexes at the infinity, the original problem mentioned at the beginning of this introduction can be formulated in the following manner: \textit{To determine properties of convex bodies imposing conditions on the sections of cones where $K$ is inscribed and whose apexes are contained in a hyperplane}. Naturally, we can replace in the aforesaid problem the condition that the set of apexes is situated in a hyperplane by the condition that they are contained in a hypersurface $S$. In particular, we can assume that $S$ is the boundary of a convex body $M\subset \Rn$ such that $K\subset \inte M$. An interesting example of this type is the well known Matsuura's Theorem \cite{Matsuura}:

\textbf{Matsuura's theorem.} Let $K\subset\mathbb R^3$ be a convex body and let $S$ be a closed convex surface which contains $K$ in its interior. If the support cone of $K$ from every point in $S$ is a right circular cone then $K$ is a Euclidean ball.

We are interested in the following problem for cones.

\textbf{Problem 2.} \emph{Given a subgroup $G$ of the general linear group $\GL(\mathbb{R},n)$ and a hypersurface $S$ which is the image of an embedding of $\Sn$, to determine the convex bodies $K\subset \Rn$, $n\geq 3$, such that for every couple of different cones $\Lambda, \Gamma$, circumscribed to $K$ and with apexes in $S$, there exists an element $\Phi \in G$ such that $\Phi(\Lambda)= \Gamma$.}

A particularly interesting case of Problem 2 is when $G$ is equal to  $O(\R,n)$, i.e., \textit{we know that all the cones which circumscribes $K$ and with apexes in $S$ are congruent}. Very recently, S. Myroshnychenko \cite{Myros} has proved the following related result: let $P$ and $Q$ be polytopes contained in the interior of the ball $B_r(n)$, $n\geq 3$, and assume that from every point in the sphere $r\mathbb S^{n-1}$ the support cones of $P$ and $Q$ are congruent, then $P=Q$.

We denote by $T(\R,n)$ the family of translations in $\Rn$. The main result of this work was inspired by Problem 2, however, we involve $T(\R,n)$ which is not a subgroup of $GL(\R,n)$, nevertheless it is an isometry of $\Rn$. Our main theorem claims that if $K \subset \mathbb{R}^n$, $n\geq 3$, is a strictly convex body and $L$ is a hypersurface, which is the image of an embedding of the sphere $\mathbb{S}^{n-1}$, $K \subset \inte L$, and for every $x\in L$, there exists $y\in L$ and $\Phi \in T(\R,n)$ such that $C_y=\Phi(C_y)$, then $K$ and $L$ are centrally symmetric and concentric.

More precisely, we are going to prove the following theorem.
  
\begin{theorem}\label{rosa}
Let $K\subset \Rn,$ $n\geq 3,$ be a strictly convex body and let $L$ be hypersurface which is an embedding of $\Sn$ such that $K\subset \inte L$. Suppose that for every $x\in L$ there exist two points $y\in L$ and $p\in \Rn$ such that
\begin{eqnarray}\label{apretada}
C_y=p+C_x.
\end{eqnarray}
Then $K$ and $L$ are centrally symmetric and concentric.
\end{theorem}

\section{Proofs and auxiliary results}
Let $\mathbb{R}^{n}$ be the Euclidean space of dimension $n$ endowed with the usual interior product 
$\langle \cdot, \cdot\rangle : \mathbb{R}^{n} \times \mathbb{R}^{n} \rightarrow \R$. We take a orthogonal coordinate system $(x_1,...,x_{n})$ for  $\mathbb{R}^{n}$. Let  
$B_r(n)=\{x\in \mathbb{R}^{n}: ||x||\leq r\}$ be the $n$-ball of radius $r$ centered  at the origin, and let $r\mathbb{S}^{n-1}=\{x\in \mathbb{R}^{n}: ||x|| = r\}$ be its boundary. For $u \in \mathbb{S}^{n-1}$ and a non-negative $s\in \R$, we denote by $\Pi(u,s)$ the hyperplane $\{x\in \mathbb{R}^{n} | \langle u, x\rangle = s \}$ whose unit normal vector is $u$ and by $\Pi^*(u,s)$ the open half-space 
$\{x\in \mathbb{R}^{n} | \langle u, x\rangle <s \}$. In particular, $\Pi(x,0)$ is denoted by $x^{\perp}$. For the points $x,y \in \Rn$ we will denote by $\aff\{x,y\}$ the affine hull of $x$ and $y$ and by $[x,y]$ the line segment with endpoints $x$ and $y$. A \textit{convex hypersurface} is the boundary of a convex body in $\Rn$. As usual $\inte K$, $\bd K$ will denote the interior and the boundary of the convex body $K$. An \textit{embedding} of $\mathbb{S}^{n-1}$ in $\Rn$ is a map $\alpha: \mathbb{S}^{n-1} \rightarrow \Rn$ such that $\alpha$ is homeomorphic onto its image.  

Let $M \subset \mathbb R^n$ be a set and let $q \in \mathbb R^n$ be a point. The set $M' = 2q-M$ is said to be centrally symmetric to $M$ with respect to the center $q$. Let $S:\Rn \rightarrow \Rn$ be the map such that $x\mapsto -x$. The following is en easy consequence of the definition of centrally symmetric sets.

\begin{lemma}\label{cosota}
Let $\Lambda$ be the image of an embedding  $\beta: \mathbb{S}^{n-2}\rightarrow \Rn$. Then $S(\Lambda)$ is a translate of $\Lambda$ if and only if $ \Lambda$ is centrally symmetric.
\end{lemma}

\begin{lemma}\label{dorron}
For every $x\in L$, the set $\Lambda=C_x\cap C_y$ is an embedding of $\mathbb{S}^{n-2}$, where $y$ satisfies (\ref{apretada}). On the other hand, the set $\Lambda$ is centrally symmetric and, if we choose the center of 
$\Lambda$ as the origin of coordinates, then 
\begin{eqnarray}\label{salsa}
-C_x=C_y.
\end{eqnarray}
\end{lemma}

\emph{Proof.}  Is easy to see that for a strictly convex body $K\subset \Rn$ and a point $x\in \Rn \setminus K$ the set $C_x\cap \bd K$, the \textit{graze} of $K$ with respect to $x$, is the image of an embedding of $\mathbb{S}^{n-2}$. On the other side, we can define a map $\phi:C_x\cap \bd K \rightarrow \Lambda$ such that $\phi$ is injective, one-to-one and bi-continuous, i.e., we can find an homeomorphism between $C_x\cap \bd K$ and 
$\Lambda$. Thus  $\Lambda$ is the image of an embedding of $\mathbb{S}^{n-2}$.

Now, note that every double-cone is a centrally symmetric set with respect to its apex. Since there is $p\in\mathbb R^n$ such that $C_y=p+C_x$, and $C_y$ is a symmetric image of $C_x$ with respect to the midpoint of $[x,y]$, we have that $\Lambda =C_x\cap C_y$ is a centrally symmetric set with center at $\frac{x+y}{2}.$ If we choose the origin of coordinates at $\frac{x+y}{2}$ then $C_y=-C_x$. \fin    

\begin{lemma}\label{punsada}
Let $x,y \in L$ such that (\ref{apretada}) holds. Let $\Pi_1,$ $\Pi_2$, $\Gamma_1,$ $\Gamma_2$, be supporting hyperplanes of $K$ such that $x\in \Pi_1\cap \Pi_2$, $\Gamma_i$ is parallel to $\Pi_i,$  $i=1,2,$ and $\Pi_1\cap \Pi_2$ is not a supporting $(n-2)$-dimensional plane of $L$. Denote by $k_i$ the contact point of $\Pi_i$ with $\bd K, i=1,2$ and by $l_i$ the contact point of $\Gamma_i$ with $\bd K, i=1,2$. Then the segments $[k_1,k_2]$ and $[l_1,l_2]$ are parallel.
\end{lemma}

\begin{figure}[H]
    \centering
    \includegraphics[width=0.58\textwidth]{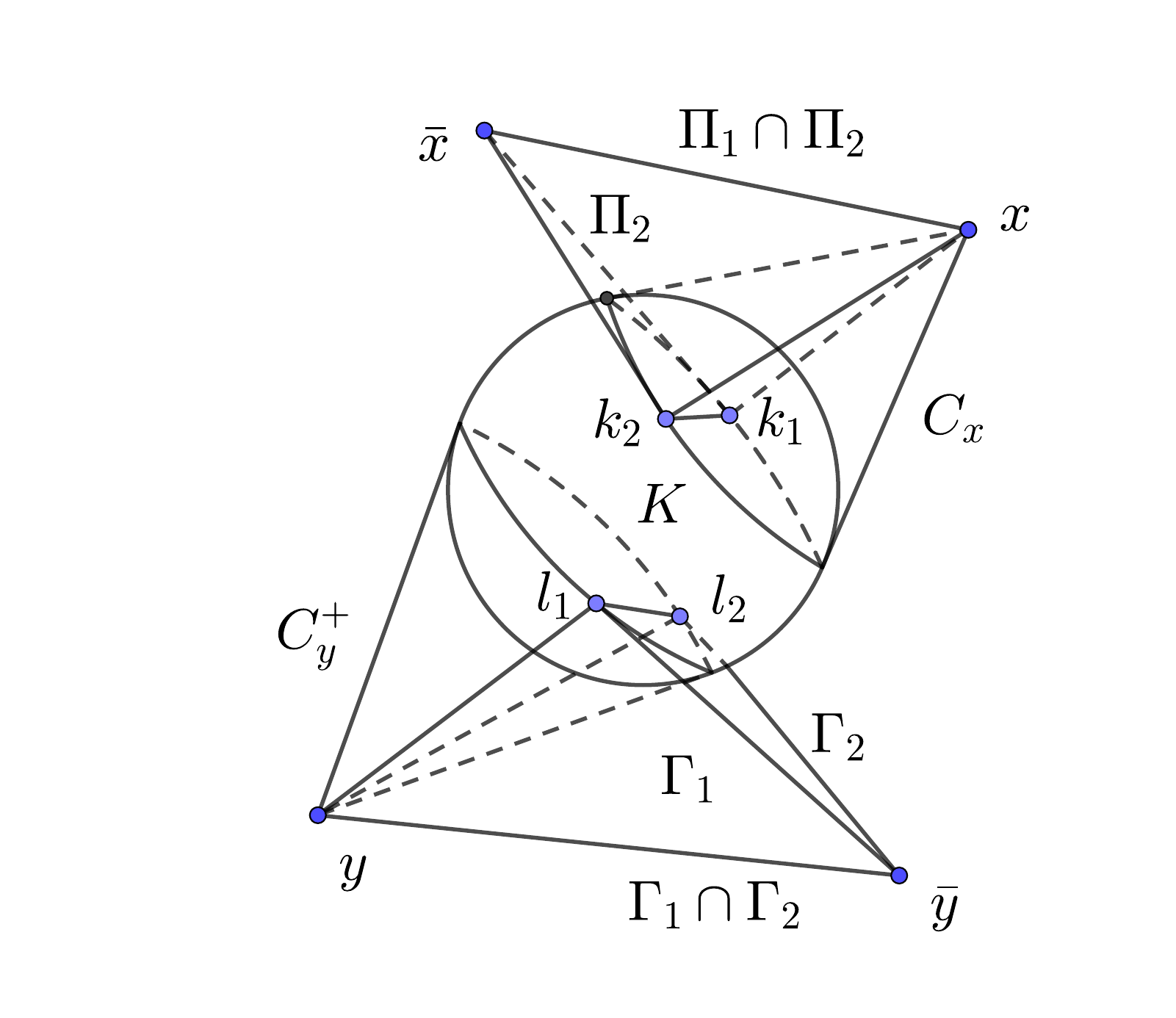}
    \caption{}
    \label{suave}
\end{figure}

\emph{Proof.} Let $\bar{x}$ be any point in $\Pi_1\cap \Pi_2\cap L$, with $\bar{x}\neq x$. In virtue of the hypothesis there exists 
$\bar{y}\in L$ and $q\in \mathbb R^n$ such that
\begin{eqnarray}\label{cuartito}
C_{\bar{y}}=q+C_{\bar{x}}.
\end{eqnarray}
By Lemma \ref{dorron} and relation (\ref{salsa}), we conclude that $y\in \Gamma_1\cap \Gamma_2$ and by (\ref{apretada}) it follows that
\[
\aff\{y,l_1\}=p+\aff\{x,k_1\} \textrm{ }\textrm{ and  }\textrm{ }\aff\{y,l_2\}=p+\aff\{x,k_2\},
\]
i.e.,
\begin{eqnarray}\label{vivaldi}
\aff\{y,l_1,l_2\}=p+\aff\{x,k_1,k_2\} 
\end{eqnarray}
(See Fig. \ref{suave}). Analogously, since $\bar{y}\in \Gamma_1\cap \Gamma_2$ by Lemma \ref{dorron} and by (\ref{cuartito}) we have that
\begin{eqnarray}\label{bach}
\aff\{\bar{y},l_1,l_2\}=q+\aff\{\bar{x},k_1,k_2\}. 
\end{eqnarray}
By (\ref{vivaldi}) and (\ref{bach}) it follows
\[
\aff\{l_1,l_2\}=\aff\{y,l_1,l_2\} \cap \aff\{\bar{y},l_1,l_2\}=
\]
\[
=[q+\aff\{x,k_1,k_2\} ]\cap [q+\aff\{\bar{x},k_1,k_2\}]=
\]
\[
=q+[\aff\{x,k_1,k_2\} \cap \aff\{\bar{x},k_1,k_2\}]
\]
\[
=q+\aff\{k_1,k_2\}.
\] \fin

\begin{lemma}\label{rojo}
Let $C\subset \Rn$ be a convex cone with apex at the origin $O$ and let $\Delta_1, \Delta_2: \mathbb{S}^{n-2}\rightarrow \Rn$ be two embeddings of $\mathbb{S}^{n-2}$ in $\bd C$. Suppose that for every $u\in \Sn$ such that $\lin \{u\}\subset \bd C$, the sets $\lin \{u\}\cap \Delta_1$ and $\lin \{u\}\cap \Delta_2$ consist of only one point; we denote by $\alpha(u)$  and by $\beta(u)$, respectively, such intersections. 
Suppose that for every two such $u,v\in \Sn, u\not=v,$ the line segment $[\alpha(u),\alpha(v)]$ is parallel to the segment $[\beta(u),\beta(v)]$. Then $\Delta_1,\Delta_2$ satisfies the relation
\begin{eqnarray}\label{gabi} 
\Delta_2=H(\Delta_1),
\end{eqnarray}
where $H:\Rn \rightarrow \Rn$ is an homothety with center of homothety at $O$.
\end{lemma}

\emph{Proof.} Let $H:\Rn \rightarrow \Rn$ be an homothety with center at $O$ such that $\Delta_2 \cap H(\Delta_1)\not=\emptyset$. Let $x\in \Delta_2 \cap H(\Delta_1)$ and let $y\in \Delta_2$ such that $y\not=x$. In virtue of the first part of hypothesis $\lin\{y\}\not=\lin\{x\}$. We denote by $z$ the point $H(\Delta_1) \cap \lin\{y\}$. It follows that $z\not=x$, otherwise $\lin\{y\}=\lin\{x\}$. We denote by $u,v$ the unit vectors $\frac{x}{\|x\|},\frac{y}{\|y\|}$, respectively. Then
\[
x=H(\alpha(u)), z=H(\alpha(v)), x=\beta(u), y=\beta(v). 
\]
In virtue of the hypothesis, the line segment $[\alpha(u),\alpha(v)]$ is parallel to the segment $[\beta(u),\beta(v)]$. Hence the line segment $[H(\alpha(u)),H(\alpha(v))]$ is parallel to the segment $[\beta(u),\beta(v)]$  and then $y=z$. Therefore, the relation (\ref{gabi}) holds. \fin

\begin{lemma}\label{morado}
Let $x$ and $y$ be points in $L$ such that (\ref{apretada}) holds. We denote by 
$\Omega_x$ and $\Omega_y$ the intersections $C_x\cap \bd K$ and $C_y\cap \bd K$, respectively. Then $\Omega_x$ and $\Omega_y$ are inversely homothetic. 
\end{lemma}

\emph{Proof.} In virtue of the hypothesis $\Delta_2=p+\Omega_x\subset C_y$. Denote by $\Delta_1$ the reflection of $\Omega_y$ with respect to $y$. Consider any two points $a,b\in\Omega_x$ and let $\Pi_a$ and $\Pi_b$ be the support hyperplanes of $K$ through $a$ and $b$, respectively, which by hypothesis pass through the point $x$. In other words, $\Pi_a$ and $\Pi_b$ are support hyperplanes of the cone $C_x$. Suppose the hyperplanes $\Pi_{a'}$ and $\Pi_{b'}$, parallel to $\Pi_a$ and $\Pi_b$ through the point $y$, touch $K$ at the points $a'$ and $b'$, respectively. Since the cone $C_y$ is a translate of the cone $C_x$, we have that $\Pi_{a'}$ and $\Pi_{b'}$ are support hyperplanes of the cone $C_y$. Without loss of generality, we suppose the origin is at the point $y$ and consider the unitary vectors $u$ and $v$ in direction of the lines $a'y$ and $b'y$, respectively. Consider the points $\alpha(u),$ $\alpha(v)$, $\beta(u)$, and $\beta(v)$, as defined in Lemma \ref{rojo}. By Lemma \ref{punsada} we have that $[a,b]$ and $[a',b']$ are parallel and so $[\alpha(u),\alpha(v)]$ is parallel to $[\beta(u),\beta(v)]$. We have the conditions of Lemma \ref{rojo}, hence $\Delta_1$ is directly homothetic to $\Delta_2$, therefore, $\Omega_x$ is inversely homothetic to $\Omega_y$. \fin

\begin{figure}[H]
    \centering
    \includegraphics[width=1.1\textwidth]{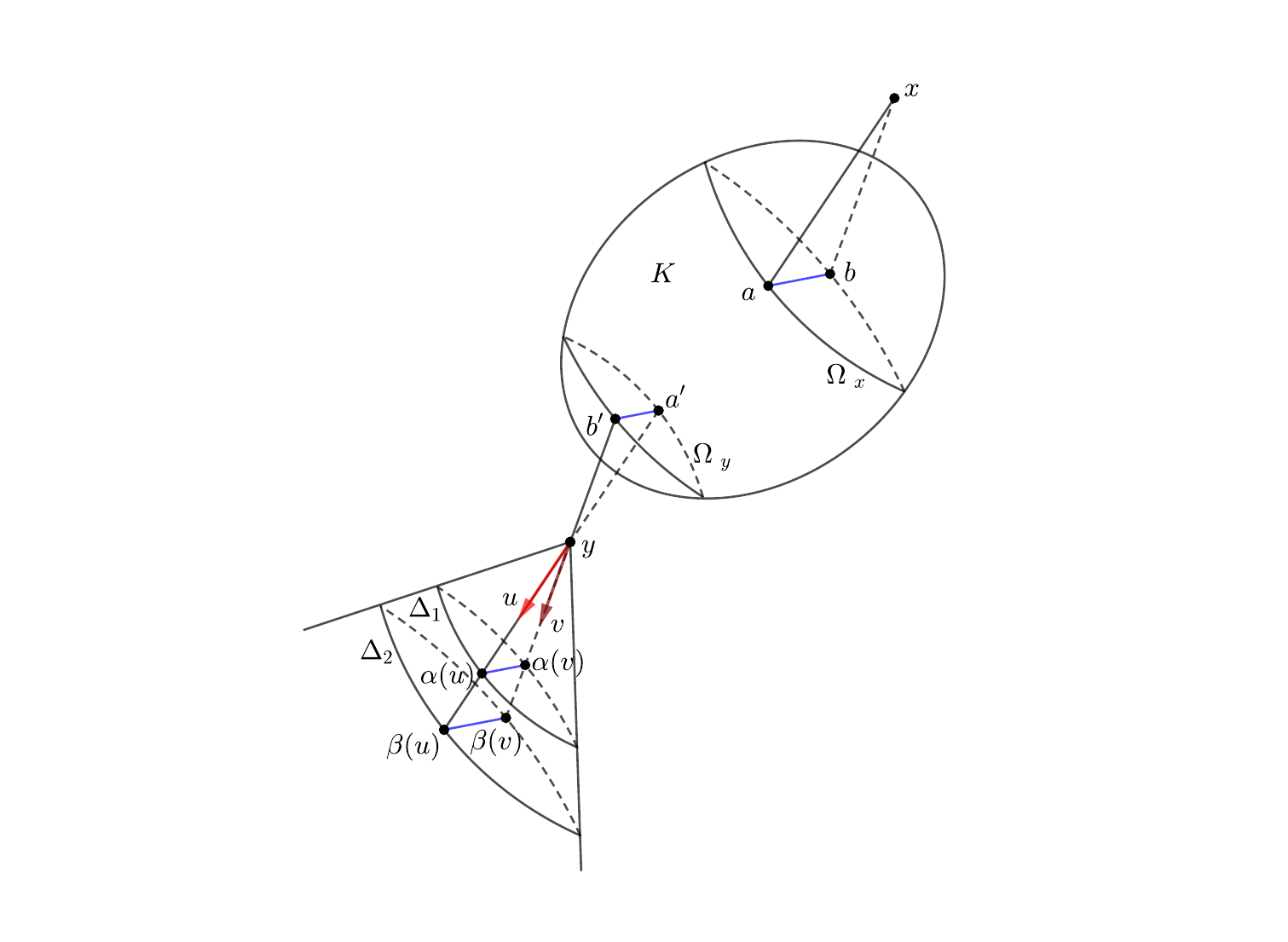}
    \caption{}
    \label{tamal}
\end{figure}

\textbf{Remark 1.} From this lemma we can deduce that every pair of parallel support hyperplanes of $K$ cut in $L$ sections which are inversely homothetic, and moreover, the center of homothety is in the segment that joins the points of contact between $K$ and these hyperplanes. With this condition and condition (\ref{apretada}), we could continue with the proof that $K$ and $L$ are centrally symmetric and concentric, however, we will continue with the proof in a little different direction. Nevertheless, this suggest that the following conjecture have chance to be true.

\textbf{Conjecture 1.} Let $K$ and $L$ be convex bodies in $\mathbb R^n$, $n\geq 3$, with $K\subset \text{int} L$. Suppose every pair of parallel support hyperplanes of $K$ cut in $L$ sections which are inversely homothetic with center of homothety in the segment that joins the points of contact between $K$ and these hyperplanes. Then $K$ and $L$ are centrally symmetric and concentric.

\begin{lemma}\label{gansito} 
Let $x,y\in L$ be two points such that (\ref{apretada}) holds. Consider the center of $\Lambda=C_x\cap C_y$ is at the origin. If $\Omega_x$ and $\Omega_y$ have ratio of homothety equal to $-1$, then
\begin{eqnarray}\label{jazz}
-\Omega_x=\Omega_y.
\end{eqnarray}
\end{lemma}

\emph{Proof.} By (\ref{salsa}) we have that $-\Omega_x\subset -C_x=C_y$. By the hypothesis 
$-\Omega_x$ and $\Omega_y$ are either translated or equal. However, Since $-\Omega_x, \Omega_y\subset C_y$ they can not be translated. Thus $-\Omega_x=\Omega_y$. \fin

\begin{lemma}\label{gelatina} 
Let $x,y,\bar{x},\bar{y}\in L$ be points such that (\ref{apretada}) holds for each pair $x,y$ and $\bar{x},\bar{y}$. Suppose that $\Omega_x \cap \Omega_{\bar{x}}$ is non-empty, it has more than one point and the ratio of homothety between $\Omega_x$ and $\Omega_y$ is -1. Then the ratio of homothety between $\Omega_{\bar{x}}$ and $\Omega_{\bar{y}}$ is -1 and the center $c$ of $\Lambda$ coincides with the center $\bar{c}$ of $\bar {\Lambda}$.
\end{lemma} 

\emph{Proof.} Following the notation of Lemma \ref{punsada}, let $\Pi_1,$ $\Pi_2$, $\Gamma_1,$ and $\Gamma_2$ be supporting hyperplanes of $K$ such that $x,\bar{x}\in \Pi_1\cap \Pi_2$ and $\Gamma_i$ is parallel to $\Pi_i,$ for $i=1,2$. We denote by $k_i$ the contact point of $\Pi_i$ with $\bd K,$ for $i=1,2$ and by $l_i$ the contact point of $\Gamma_i$ with $\bd K,$ for $i=1,2$. By Lemma \ref{punsada} we have that the segments $[k_1,k_2]$ and $[l_1,l_2]$ are parallel. On the other side, since $k_1,k_2\in \Omega_x$ and $l_1,l_2\in \Omega_y$ and the ratio of homothety between $\Omega_x$ and $\Omega_y$ is -1, which is indeed the ratio of homothety between the segments $[k_1,k_2]$ and  $[l_1,l_2]$ such that $k_i$ corresponds to $l_i$ (for $i=1,2$), we have that $\|k_1-k_2\|=\|l_1-l_2\|$. The center of the inverse homothety between $\Omega_x$ and $\Omega_y$ is the center of the parallelogram $k_2k_1l_1l_2,$ which is precisely $c=y-p/2.$ We also have that the center of symmetry of $\Lambda$ is $c$. Now, since $k_1,k_2\in \Omega_{\bar{x}}$ and $l_1,l_2\in \Omega_{\bar{y}}$, it follows that the ratio of homothety between $\Omega_{\bar{x}}$ and $\Omega_{\bar{y}}$ is $-1$ with center at $c$ as well. \fin

Let $N\subset L$ be the set of points $x\in L$ such that, if the point $y\in L$ has the property that for $x,y$  (\ref{apretada}) holds, then the curves $\Omega_x$ and $\Omega_y$ are inversely homothetic with factor of homothety equal to $-1$. We will prove that indeed $N=L$.

\begin{lemma}\label{perrita} 
We have that $N=L$. 
\end{lemma}

\emph{Proof.} Consider a pair of points $u,v\in L$ such that (\ref{apretada}) holds for $u,v$ and $u,v\in N$. The existence of this pair of points is guaranteed by an standard argument of continuity, however, we give some details of this argument for the interested reader in the Appendix. 

We denote by $D_u$ the set of points $z\in L$ such that $[u,z]\cap K\neq\emptyset$ and denote by $L_u$ the set $L\setminus D_u$. Let $a$ be any point in $L_u$, with $a\neq u$. Since $\Omega_a\cap\Omega_u\neq\emptyset$, by Lemma \ref{gelatina} we have that $a\in N$. Associated with $a$ we have the open region $L_a$, and as before, every point $z\in L_a$ is also in $N$. Continuing with this procedure we construct an open cover of $L$ such that every open set in this cover belong to $N$. Since $L$ is a compact set, we may consider indeed a finite subcover of $L$, hence all the points in $L$ belong to $N$. \fin
  
\emph{Proof of Theorem \ref{rosa}}. First, by Lemma \ref{gelatina} we have that all the curves $\Lambda$ are concentric at $c$. From here we get that the support function $h_k:\mathbb{S}^{n-1} \rightarrow \mathbb{R}$ of $K$ is symmetric, i.e., $h_k(u)=h_k(-u)$, that is, $K$ is centrally symmetric with center at $c$.

Finally, we prove that $L$ is centrally symmetric and concentric with $K$. We choose a system of coordinates such that the origin $O$ is the centre of $K$. We observe that if a convex body $M\subset \Rn$ is centrally symmetric with centre at $O$ and $x\in \Rn \setminus M$, there exist $y\in \Rn \setminus K$ and $p\in \Rn$, both unique, such that the relation
\[
C_y=p+C_x
\]
holds, it is clear that $y=-x$. If $L$ were not be centrally symmetric and concentric with $K$, then would exist $x_0\in L$ such that $-x_0\not\in L$. Thus would not exist $y_0\in L$ and $p\in \Rn$ such that $C_{y_0}=p+C_{x_0}$ which would contradict the hypothesis of Theorem \ref{rosa}. Therefore, $L$ is centrally symmetric with center at $O$. \fin

\textbf{Remark 2.} The main theorem proved in this article is not true in the plane, even if all the support cones are congruent. To see this is sufficient to consider a convex body $K\subset\mathbb R^2$ and its \emph{isoptic curve} $K_{\alpha}$, i.e., the set of points $z\in\mathbb R^2$ such that $K$ is seen under the constant angle $\alpha$ from $z$ (see for instance \cite{Cies_Mier_Moz}, \cite{Green}).

\section{Appendix}
Here we prove that there exists two points $u,v\in L$ such that (\ref{apretada}) holds for $u,v$ and $u,v\in N$. 

\begin{figure}[H]
    \centering
    \includegraphics[width=.7\textwidth]{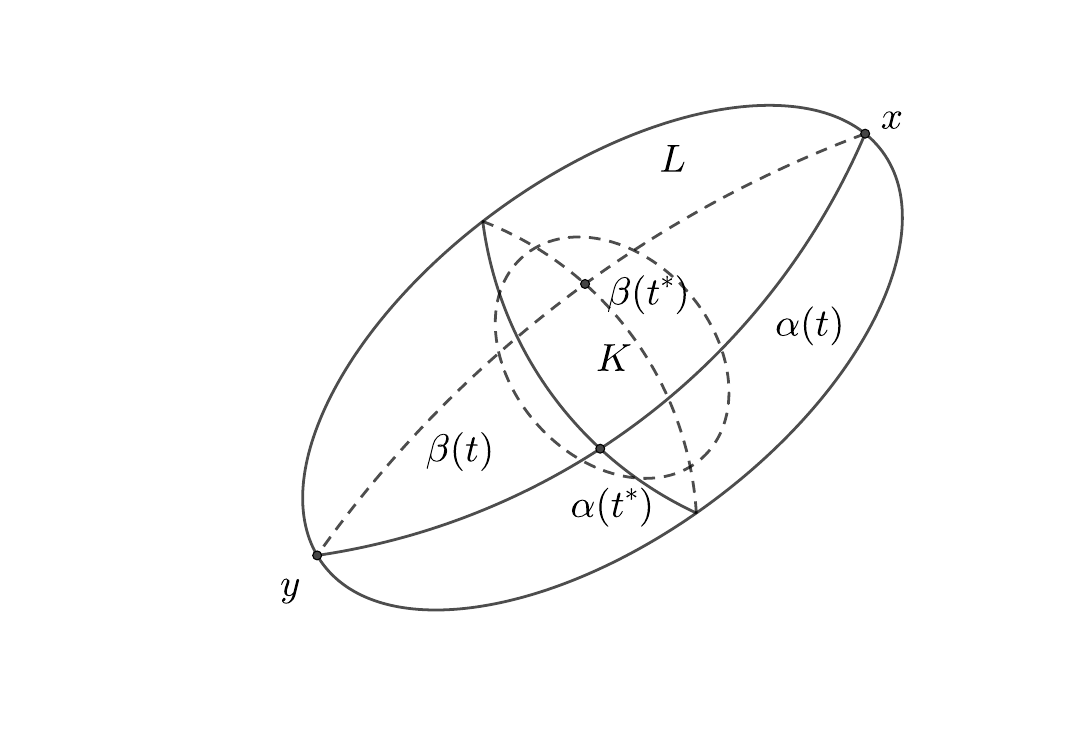}
    \caption{}
    \label{guera}
\end{figure}

Let $\alpha:[0,1] \rightarrow L$ be a continuous curve such that $\alpha(0)=x$ and $\alpha(1)=y$. In virtue of the hypothesis there exists a continuous curve $\beta:[0,1] \rightarrow L$ such that $\beta(0)=y$ and $\beta(1)=x$ and, for every $t\in[0,1]$, the points $\alpha(t)$ and $\beta(t)$ are such the relation (\ref{apretada}) holds for them, i.e., 
\begin{eqnarray}\label{blues}
C_{\beta(t)}=p(t)+ C_{\alpha(t)},
\end{eqnarray}
where $p(t)\in \Rt$. We  define $d(t)$ as the number $ \diam \Omega_{\beta(t)}-\diam \Omega_{\alpha(t)}$. By Lemma \ref{morado}, 
\begin{eqnarray}\label{todo}
\diam \Omega_{\alpha(t)}=r(t)\diam \Omega_{\beta(t)}. 
\end{eqnarray}
Thus 
\[
d(t)=(1-r(t))\diam \Omega_{\beta(t)}.
\]
Since the surfaces $L$ and $\bd K$ are continuous and $L\cap K=\emptyset$, it follows that the functions $r,d:[0,1]\rightarrow \R^+$ are continuous. By our assumption 
$\diam \Omega_{\alpha(0)}\not=\diam \Omega_{\beta(0)}$, say 
\begin{eqnarray}\label{recarga}
\diam \Omega_{\alpha(0)}<\diam \Omega_{\beta(0)},
\end{eqnarray}
replacing (\ref{todo}), for $t=0$, in (\ref{recarga}) we get $r(0)<1$, that is, $d(0)>0$. On the other hand,
in virtue that $\alpha(0)=x=\beta(1)$ and $\alpha(1)=y=\beta(0)$, from (\ref{recarga}) it follows 
\[
\diam \Omega_{\beta(1)}<\diam \Omega_{\alpha(1)},
\]
i.e., $1<r(1)$ and $d(1)<0$. Consequently, by the mean value theorem, there exists $t^*\in [0,1]$ such that $d(t^*)=0$, that is, $r(t^*)=1$. Hence 
\begin{eqnarray}\label{blanquita}
\diam \Omega_{\alpha(t^*)}=\diam \Omega_{\beta(t^*)}.
\end{eqnarray} 
(See Fig. \ref{guera}). Thus $u=\alpha(t^*)$ and $v=\beta(t^*)$ are the point what we are looking for.

\end{document}